\documentclass[12pt,leqno]{article}
\usepackage[centertags]{amsmath}
\usepackage{amsthm}
\usepackage{amssymb}
\usepackage{mathtools}
\numberwithin{equation}{section}

\newtheorem{theorem}{Theorem}
\newtheorem{lemma}{Lemma}

\newtheorem{rem}{Remark}

\newcommand{\Landau}{\ensuremath{\mathcal O}}

\newcommand{\eps}{\varepsilon}

\begin{document}

\title {Error estimates in balanced norms of finite element methods
on Shishkin meshes for reaction-diffusion problems}

\author{  Hans-G. Roos}

\date{}

\maketitle

\begin{abstract} Error estimates of finite element methods for reaction-diffusion
problems are often realized in the related energy norm. In the singularly perturbed
case, however, this norm is not adequate. A different scaling of the $H^1$ seminorm
leads to a balanced norm which reflects the layer behavior correctly. We discuss
also anisotropic problems, semilinear equations, supercloseness and a combination
technique.

\end{abstract}


\section{Introduction}

We shall examine the finite element method for the numerical
solution of the singularly perturbed linear elliptic boundary value
problem
\begin{subequations}\label{1.1}
  \begin{align}
Lu\equiv -\varepsilon \Delta u +cu&=f
\qquad \text{in } \Omega=(0,1)\times(0,1)\\
u&=0 \qquad \text{on } \partial \Omega,
   \end{align}
\end{subequations}
where $0<\varepsilon \ll 1$ is a small positive parameter,
$c>0$ is (for simplicity) a positive constant and $f$ is sufficiently smooth.

It is well-known that the problem has a unique solution $u\in V=H_0^1(\Omega)$
which satisfies the stability estimate in the related energy norm
\begin{equation}\label{1.2}
\|u\|_{\varepsilon}\coloneqq\varepsilon^{1/2}|u|_1+\|u\|_0\preceq \|f\|_0.
\end{equation}
Here we used the following notation: if $A\preceq B$, there exists a (generic)
constant $C$ independent of $\varepsilon$ (and later also of the mesh used)
such that $A\le C\,B$. Moreover for $D \subset \Omega$ we denote by $\|\cdot\|_{0,D}$, $\|\cdot\|_{\infty,D}$ and $|\cdot|_{1,D}$ the standard norms in $L_2(D)$, $L_\infty(D)$ and the standard seminorm in $H^1(D)$, respectively. We shall omit the notation of the domain in the case $D = \Omega$. Similarly, we want to use the notation $(\cdot,\cdot)_D$ for the inner product in $L_2(D)$ and abbreviate $(\cdot,\cdot)_\Omega$ to $(\cdot,\cdot)$.

Moreover, the error of a finite element approximation $u^N\in V^N\subset V$
satisfies
\begin{equation}\label{1.3}
\|u-u^N\|_{\varepsilon}\preceq \min_{v^N\in V^N}\|u-v^N\|_{\varepsilon}.
\end{equation}

When linear or bilinear elements are used on a Shishkin mesh (see
Section 2), one can prove under certain additional assumptions
concerning $f$ for the interpolation error of the Lagrange
interpolant $u^I\in V^N$ on Shishkin meshes
\begin{equation}\label{1.4}
\|u-u^I\|_{\varepsilon}\preceq \left(\varepsilon^{1/4}N^{-1}\ln N+N^{-2}\right)
\end{equation}
(see \cite{DR97} or \cite{RST96}). It follows that the error $u-u^N$ also satisfies such an estimate.

However, the typical boundary layer function $\exp(-x/\varepsilon^{1/2})$
measured in the norm $\|\cdot\|_{\varepsilon}$ is of order $\Landau(\varepsilon^{1/4})$.
Consequently, error estimates in this norm are less valuable as for convection
diffusion equations where the layers are of the structure $\exp(-x/\varepsilon)$.
Wherefore we  ask the fundamental question:\\
\emph{Is it possible to prove error estimates in the balanced norm}
\begin{equation}\label{1.5}
\|v\|_{b}\coloneqq\varepsilon^{1/4}|v|_1+\|v\|_0\quad ?
\end{equation}
In Section 2 we will repeat an basic idea to prove error estimates in a balanced norm
and extend the approach to semilinear problems and anisotropic equations. Supercloseness
and a combination technique are discussed in Section 3. Finally we present a direct
mixed method in Section 4.

\section{The basic error estimate in a balanced norm and some extensions
}\label{sec:2}
The mesh $\Omega^N$ used is the tensor product of two
one-dimensional piecewise uniform Shishkin meshes. I.e.,
$\Omega^N=\Omega_x\times\Omega_y$, where $\Omega_x$ (analogously
$\Omega_y$) splits $[0,1]$ into the subintervals
 $[0,\lambda_x]$, $[\lambda_x,1-\lambda_x]$ and $[1-\lambda_x,1]$. The mesh
 distributes $N/4$ points equidistantly within each of the subintervals
 $[0,\lambda_x]$, $[1-\lambda_x,1]$ and the remaining points within the third
 subinterval. For simplicity, assume
 \[
 \lambda=\lambda_x=\lambda_y=\min\{1/4, \lambda_0\sqrt{\varepsilon/c^*}\ln N\}
 \quad\text{with }\lambda_0=2 \text{ and } c^*<c.
 \]
 We use for the step sizes
 \[
 h \coloneqq \frac{4\lambda}{N}\quad{\rm and}\quad H\coloneqq\frac{2(1-2\lambda)}{N}.
 \]

Let $V^N\subset H_0^1(\Omega)$ be the space of bilinear finite elements
on $\Omega^N$ or the space of linear elements over a triangulation obtained
from $\Omega^N$ by drawing diagonals.

A standard formulation of problem \eqref{1.1} reads: find $u \in V$, such that
\begin{gather}\label{eq:standard}
	\varepsilon (\nabla u,\nabla v) + c (u,v) = (f,v)\quad \forall v \in V.
\end{gather}
By replacing $V$ in \eqref{eq:standard} with $V^N$ one obtains a standard discretization
that yields the FEM-solution $u^N$.

As we mentioned already in the Introduction, certain assumptions on $f$
allow a decomposition of $u$ into smooth components and layer terms
such that the following estimates for the interpolation error of the
Lagrange interpolant hold true (see \cite{DR97} or \cite{RST96}):

\begin{alignat}{2}
\|u-u^I\|_0 &\preceq N^{-2},&\quad \varepsilon^{1/4}|u-u^I|_1 &\preceq N^{-1}\ln N \label{2.1}
\shortintertext{and}
\|u-u^I\|_{\infty, \Omega_0} &\preceq N^{-2},&\quad \|u-u^I\|_{\infty,\Omega\setminus \Omega_0} &\preceq (N^{-1}\ln N)^2, \label{2.2}
\end{alignat}
here $\Omega_0=(\lambda_x,1-\lambda_x)\times(\lambda_y,1-\lambda_y)$. Let us also introduce
$\Omega_f:=\Omega\setminus \Omega_0$.

Instead of the Lagrange interpolant we use in our error analysis the $L_2$
projection $\pi u\in V^N$ from $u$. Based on
\[u-u^N=u-\pi u+\pi u-u^N
\]
we estimate $\xi\coloneqq\pi u-u^N$:
\[
\|\xi\|_\varepsilon^2\preceq \varepsilon|\nabla
\xi|_1^2+c\,\|\xi\|_0^2=
   \varepsilon(\nabla(\pi u-u),\nabla \xi)+c\,(\pi u-u,\xi).
   \]
Because $(\pi u-u,\xi)=0$, it follows
\begin{equation}\label{2.3}
|\pi u-u^N|_1\preceq |u-\pi u|_1.
\end{equation}
If we now could prove a similar estimate as \eqref{2.1} for the error of the
$L_2$ projection, we obtain an estimate in the balanced norm because we have
already an estimate for $\|u-u_N\|_0$ in \eqref{1.4}.
\begin{lemma}
Assuming the validity of \eqref{2.1} and \eqref{2.2}, the error of
the $L_2$ projection on the Shishkin mesh satisfies
\begin{equation}\label{2.4}
\|u-\pi u\|_\infty\preceq \|u- u^I\|_\infty
,\quad \varepsilon^{1/4}|u-\pi u|_1\preceq N^{-1}(\ln N)^{3/2}.
\end{equation}
\end{lemma}
The proof uses
the $L_\infty$-stability of the $L_2$ projection on our mesh
\cite{Os08}. Inverse inequalities are used to move from estimates in
$W^1_\infty$ to $L_\infty$, for details see \cite{RS15}.

From \eqref{2.3} and Lemma 1 we get
\begin{theorem}
Assuming \eqref{2.1} and \eqref{2.2}, the error of the Galerkin finite
element method with linear or bilinear elements on a Shishkin mesh satisfies
\begin{equation}\label{2.5}
\|u-u^N\|_{b}\preceq N^{-1}(\ln N)^{3/2}+N^{-2}.
\end{equation}
\end{theorem}

Remark that for $Q_k$ elements with $k>1$ one can get an analogous result
\[
\|u-u^N\|_{b}\preceq N^{-k}(\ln N)^{k+1/2}+N^{-(k+1)}
\]
because on tensor product meshes the $L_2$ projection is as well $L_\infty$ stable
(see \cite{CT87} for the one-dimensional result on arbitrary meshes, on tensor
product meshes the statement follows immediately).

It is easy to modify the basic idea to the
singularly perturbed semilinear elliptic boundary value
problem
\begin{subequations}
  \begin{align}
Lu\equiv -\varepsilon \Delta u +g(\cdot,u)&=0
\qquad \text{in } \Omega=(0,1)\times(0,1)\\
u&=0 \qquad \text{on } \partial \Omega.
   \end{align}
\end{subequations}
 We assume
 that  $g$ is sufficiently smooth and $\partial_2 g\ge \mu>0$. Then, the so
 called reduced problem and our given problem have a unique solution.

If $\partial \Omega$ is smooth, the solution is characterized by the typical
boundary layer for linear reaction-diffusion problems, see \cite{JF96} for the
semilinear case. If corners exist, additionally corner layers arise, see \cite{KK05}
for semilinear problems in a polygonal domain. For the analysis of finite element
methods on layer-adapted meshes we need a solution decomposition (see Remark 1.27 in
Chapter 3 of \cite{RST96}), in the semilinear case sufficient conditions for the
existence of such a decomposition are not known. Therefore we just assume the existence
of a solution decomposition.

A standard weak formulation of our semilinear problem reads: find $u \in V$, such that
\begin{gather}\label{eq:standard}
	\varepsilon (\nabla u,\nabla v) +(g(\cdot,u),v)  = 0\quad \forall v \in V.
\end{gather}
By replacing $V$ in \eqref{eq:standard} with $V^N$ one obtains a standard discretization
that yields the FEM-solution $u^N$.

If $\pi u\in V^N$ is some projection of $u$, we decompose the error into
\[u-u^N=u-\pi u+\pi u-u^N
\]
and (assuming we can control the projection error) start the error analysis from
the following relation for  $\xi:=\pi u-u^N$:

\[
\varepsilon|\nabla
\xi|_1^2+\mu\,\|\xi\|_0^2 \le \varepsilon(\nabla \xi,\nabla \xi)+(g(\cdot,\pi u)-g(\cdot,u^N),\xi)=
   \varepsilon(\nabla(\pi u-u),\nabla \xi)+(g(\cdot,\pi u)-g(\cdot,u),\xi).
   \]
If we choose $\pi u$  to be the standard interpolant of $u$, the usual error
estimate in the energy norm
\begin{equation}
\|u\|_{\varepsilon}:=\varepsilon^{1/2}|u|_1+\|u\|_0
\end{equation}
   follows:
\begin{equation}\label{ener}
\|u-u^N\|_{\varepsilon}\preceq\left(\varepsilon^{1/4}N^{-1}\ln N+N^{-2}\right)
\end{equation}

But we want again to prove an error estimate in the balanced norm
\begin{equation}
\|v\|_{b}:=\varepsilon^{1/4}|v|_1+\|v\|_0\quad .
\end{equation}

Following the basic idea from \cite{RS15}, we define $\pi u$  by
\begin{equation}\label{pro}
(g(\cdot,\pi u),v)=(g(\cdot, u),v)\quad {\rm for\,\, all}\,v\in V^N.
\end{equation}
Our assumption $\partial_2 g\ge \mu>0$ tells us immediately that $\pi u$ is well
defined, moreover
\begin{equation}\label{L2}
    \|u-\pi u\|_0\preceq \inf_{v\in V^N}\|u-v\|_0.
\end{equation}
It follows from the definition of our projection
\begin{equation}\label{2.3}
|\pi u-u^N|_1\preceq |u-\pi u|_1.
\end{equation}

For the standard interpolant $u^I$ of $u$ we have
\[
\eps^{1/4}|u-u^I|_1\preceq N^{-1}\ln N.
\]

If we now could prove a similar estimate for our projection error,
 we would obtain an estimate in the balanced norm because we have
already an estimate for $\|u-u_N\|_0$ in \eqref{ener}.

\begin{lemma}
The projection defined by \eqref{pro} is $L_\infty$ stable.
\end{lemma}

{\it Proof}: The proof is based on Taylors formula
\[
F(w)-F(v)=(\int_0^1 DF(v+s(w-v))ds)(w-v).
\]
Introducing the linear operator
\[
\triangle F(v,w):=\int_0^1 DF(v+s(w-v))ds
\]
it is obvious that
\[
\|w-v\|\le \|(\triangle F(v,w)^{-1}\|\,\|F(w)-F(v)\|.
\]
Therefore, the $L_\infty$ stability of the $L_2$ projection on our mesh
\cite{Os08} implies the  $L_\infty$ stability of our generalized projection
as well.

\begin{lemma}
 The projection error of \eqref{pro} on the Shishkin mesh satisfies
\begin{equation}\label{2.4}
\|u-\pi u\|_\infty\preceq \|u- u^I\|_\infty
,\quad \varepsilon^{1/4}|u-\pi u|_1\preceq N^{-1}(\ln N)^{3/2}.
\end{equation}
\end{lemma}
The proof works analogously as in the linear case. And, consequently, we
get the same error estimate as in Theorem 1 also in the semilinear
case.

Next we consider the anisotropic problem
\begin{subequations}\label{aniso}
  \begin{align}
 -\varepsilon  u_{xx}+u_{yy} +cu&=f
\qquad \text{in } \Omega=(0,1)\times(0,1)\\
u&=0 \qquad \text{on } \partial \Omega.
   \end{align}
\end{subequations}
Now we have only boundary layers at $x=0$ and $x=1$, the layers are of elliptic type. But the
layer terms satisfy the same estimates as in the reaction-diffusion regime \cite{Li97}. Therefore,
the estimates (2.2) and \eqref{2.2} for the interpolation error on the related Shishkin mesh remain valid,
of course, now $\Omega_0=(\lambda_x,1-\lambda_x)\times(0,1)$. Therefore, defining the energy norm by
\[
\|v\|_{\varepsilon,a}\coloneqq\varepsilon^{1/2}\|u_x\|_0+\|u_y\|_0+\|u\|_0
\]
it follows for bilinear elements
\[
\|u-u^N\|_{\varepsilon,a}\preceq \left(\varepsilon^{1/4}N^{-1}\ln N+N^{-2}\right).
\]
If we want to estimate the error in the balanced norm
\[
\|v\|_{b,a}\coloneqq\varepsilon^{1/4}\|u_x\|_0+\|u_y\|_0+\|u\|_0,
\]
we start for $\xi:=\pi u-u^N$ from
\[
\eps\|\xi_x\|_0^2\le \eps((\pi u-u)_x,
\xi_x)+
   ((\pi u-u)_y, \xi_y)+c\,(\pi u-u,\xi).
\]
Now we define in the anisotropic  case the projection onto the finite element space by
\[
  ((\pi u-u)_y, \xi_y)+c\,(\pi u-u,\xi)=0 \quad \forall \xi\in V^N.
\]
Consequently it remains to estimate for that projection $\|(\pi u-u)_x\|_0$. But the projection
satisfies
\[
  \pi v=\pi^y(\pi^x\,v),
\]
where $\pi^x$ is the one-dimensional $L_2$ projection and $\pi^y$ the one-dimensional Ritz
projection (with respect to a non-singularly perturbed operator on a standard mesh), compare
\cite{FR14}. Consequently, the projection is $L_\infty$ stable and we can repeat our
basic idea to prove estimates in the balanced norm.

Remark that in \cite{MX15} the authors use a different technique to derive estimates in the
$H^1$ seminorm for the generalized $L_2$ projection on a layer-adapted mesh.

\section{Supercloseness and a combination technique}
We come back to the linear reaction-diffusion problem
\begin{subequations}\label{1.1}
  \begin{align}
Lu\equiv -\varepsilon \Delta u +cu&=f
\qquad \text{in } \Omega=(0,1)\times(0,1)\\
u&=0 \qquad \text{on } \partial \Omega.
   \end{align}
\end{subequations}
For bilinear elements on the corresponding Shishkin mesh it is well known that we have the
supercloseness property (assuming $\lambda_0\ge 2.5$)
\begin{equation}\label{1.4}
\|u^N-u^I\|_{\varepsilon}\preceq \left(\varepsilon^{1/2}(N^{-1}\ln N)^2+N^{-2}\right).
\end{equation}
Now we ask: does there exist some projection onto the finite element space such that a supercloseness
property holds with respect to the balanced norm?

With $v_N:=u^N-\Pi u$ we start from
\[
 \varepsilon|v_N|_1^2+c\,\|v_N\|_0^2\preceq
   \varepsilon(\nabla(u-\Pi u),\nabla v_N)+c\,(u-\Pi u,v_N).
   \]
Next we use the decomposition $u=S+E$, decompose also  $\Pi u=\Pi S+\Pi E$ and use different projections
into our bilinear finite element space for $S$ and $E$. We choose:
\begin{itemize}
\item $\Pi S\in V^N$ satisfies
   \[
    (\Pi S,v)=(S,v)\quad \forall v\in V^N_0
   \]
with given values in the grid points on the boundary.
\item $\Pi E$ is zero in $\Omega_0$ and the standard bilinear interpolation operator in the fine subdomain
      with exception of one strip of the width of the fine stepsize in the transition region
      (and, of course, bilinear in that strip and globally continuous)
\end{itemize}
With this choice we obtain
\[
 \varepsilon|v_N|_1^2+c\,\|v_N\|_0^2\preceq
   \varepsilon(\nabla(u-\Pi u),\nabla v_N)+c\,(E-\Pi E,v_N)_{\Omega_f}.
   \]
In the second term we hope to get some extra power of $\eps$, in the first term we want to apply
superconvergence techniques for the estimation of the expression $(\nabla(E-\Pi E),\nabla v_N$).
First let us remark that $\Pi E$ satisfies the same estimates as the bilinear interpolant $E^I$ on
$\Omega_f$:
\[
   \|E-\Pi E\|_{0,\Omega_f}\preceq \eps^{1/4}(N^{-1}\ln N)^2
\]
and (based on Lin identities)
\[
 \varepsilon|(\nabla(E-\Pi E),\nabla v_N)|\preceq N^{-2}\eps^{3/4}|v_N|_1.
\]
It is only a technical question to prove that for our modified interpolant using the fact that $E$
is on that strip is as small as we want and that the measure of the strip is small as well.

Consequently we get
\[
|v_N|_1^2 \preceq |S-\Pi S|_1^2 +\eps^{-1/2}(N^{-1}\ln N)^4.
\]
For the $L_2$ projection of $S$ we have $\|S-\Pi S\|_\infty\preceq N^{-2}$ and
$\|S-\Pi S\|_{\infty,\Omega_f}\preceq (\eps^{1/2}N^{-1}\ln N)^2$. It follows
\[
|S-\Pi S|_{1,\Omega_0}\preceq N^{-1}, \quad |S-\Pi S|_{1,\Omega_f}\preceq \eps^{1/2}N^{-1}\ln N.
\]
Summarizing we get the supercloseness result
\[
\eps^{1/4}|u^N-\Pi u|_1\preceq \eps^{1/4}N^{-1}+(N^{-1}\ln N)^2.
\]
It is no problem to estimate the $L_2$ error.

Next we present an application of the supercloseness result to the combination technique.
We analyse the version of the combination technique presented in \cite{FLR09}, for a different
version see \cite{LMS09}. Remark that in \cite{MR15} the authors observe numerically a nice
behaviour of a combination technique in the balanced norm.

Writing $N$ for the maximum number of mesh intervals in each coordinate direction, our combination
technique simply adds or subtracts solutions that have been computed by the Galerkin FEM on
$N\times \sqrt{N}$, $ \sqrt{N}\times N$ and $ \sqrt{N}\times \sqrt{N}$ meshes. We obtain
the same accuracy as on a $N\times N$ mesh with less degrees of freedom. In the following
we use the notation of \cite{FLR09}.

In the combination technique for bilinear elements we compute a two-scale finite element
approximation $u_{\hat N,\hat N}^N$ by
\[
u_{\hat N,\hat N}^N:=u_{ N,\hat N}^N+u_{\hat N, N}^N-u_{\hat N,\hat N}^N.
\]
Later we will choose $\hat N=\sqrt{N}$. We proved (in our new notation)
\begin{equation}
\|u-u_{NN}\|_{b}\preceq N^{-1}(\ln N)^{3/2}+N^{-2}.
\end{equation}
The question is wether or not $u_{\hat N,\hat N}^N$ satisfies a similar estimate (in the
case $\hat N=\sqrt{N}$).

Analogously to $u_{\hat N,\hat N}^N$ we define $I_{\hat N,\hat N}^N E$ and $\Pi_{\hat N,\hat N}^N S$.
Then we can as follows decompose the error to estimate:
\[
u_{\hat N,\hat N}^N-u_{NN}=T_{cl,1}(S)+(\Pi_{\hat N,\hat N}^N S-\Pi_{N,N}S)+
     T_{cl,s}(E)+(I_{\hat N,\hat N}^N E-I_{N,N}E)
\]
Thus we have two terms representing the error for two-scale projection operators (related to $L_2$
projection and interpolation, respectively) and two terms which can be estimated based on our
supercloseness result:
\[
T_{cl,1}(S):=
(S_{N,\hat N}-\Pi_{N,\hat N} S)+(S_{\hat N,N}-\Pi_{\hat N,N} S)-(S_{\hat N,\hat N}-\Pi_{\hat N,\hat N} S),
\]
analogously
\[
T_{cl,2}(E):=
(E_{N,\hat N}-I_{N,\hat N} E)+(E_{\hat N,N}-I_{\hat N,N} E)-(E_{\hat N,\hat N}-I_{\hat N,\hat N} E),
\]
For the two-scale interpolation error $(I_{\hat N,\hat N}^N E-I_{N,N}E)$ the results of \cite{FLR09} remain
valid (Lemma 2.3 and 2.4, modified for the reaction-diffusion problem). For the two-scale projection error
an estimate in $L_2$ and $L_\infty$ is easy. The estimate in the seminorm $|\cdot|$ as in Section 2 follows
from an inverse inequality, applied separately in $\Omega_0$ and $\Omega_f$
. Finally we get for $\hat N=\sqrt{N}$ the estimate
\begin{equation}
\|u_{\hat N,\hat N}^N-u_{NN}\|_{b}\preceq \eps^{1/4}N^{-1/2}+N^{-1}\ln N.
\end{equation}
That means so far we can only proof the desired estimate for the combination technique if $\eps\preceq N^{-2}$.

\section{A direct mixed method}
The first balanced error estimate was presented by Lin and Stynes \cite{LS12} using a first order
system least squares (FOSLS) mixed method. For the variables $(u,\bar q)$ with $-\bar q=\nabla u$
and its discretizations on a Shishkin mesh they proved
\begin{equation}\label{est1}
\varepsilon^{1/4}|\bar q-\bar q^N|_1+\|u-u^N\|_0\preceq N^{-1}\ln N
\end{equation}
(see also \cite{ALM15} for a modified version of the method)

We shall proof that the estimate \eqref{est1} is also valid for a direct mixed
method (instead the more complicated least-squares approach from \cite{LS12}).
 Remark that Li and Wheeler \cite{LW00} analyzed the method
in the energy norm on so called A-meshes, simpler to analyze than S-meshes.

Introducing $\bar q=-\nabla u$, a weak formulation of \eqref{1.1} reads:\\
Find $(u,\bar q)\in V\times W$ such that
\begin{subequations}\label{2.1}
  \begin{align}
\eps(div\, \bar q,w)+c(u,w)&=(f,w)
\qquad \text{for all } w\in W,\\
\eps(\bar q,\bar v)-\eps(div\, \bar v,u)&=0 \qquad \text{for all } \bar v \in V,
   \end{align}
\end{subequations}
with $V=H(div,\Omega)$, $W=L^2(\Omega)$.

For the discretization on a standard rectangular Shishkin mesh (see \cite{LS12}, page 2735)
we use $(u^N,\bar q^N)\in V^N\times W^N $. Here $W^N$ is the space of piecewise constants on our rectangular mesh
and $V^N$ the lowest order Raviart-Thomas space $RT_0$. That means, on each mesh rectangle elements
of $RT_0$ are vectors of the form
\[
  (span(1,x), span(1,y)^T.
\]
Our discrete problem reads: Find $(u^N,\bar q^N)\in V^N\times W^N$ such that
\begin{subequations}\label{n2.2}
  \begin{align}
\eps(div\, \bar q^N,w)+c(u^N,w)&=(f,w)
\qquad \text{for all } w\in W^N,\\
\eps(\bar q^N,\bar v)-\eps(div\, \bar v,u^N)&=0 \qquad \text{for all } \bar v \in V^N.
   \end{align}
\end{subequations}
Setting $w:=u^N$, $\bar v:=\bar q^N$ results in the stability estimate
\begin{equation}\label{2.3}
\eps\|\bar q^N\|_0^2+\frac{c}{2}\|u^N\|_0^2\preceq \|f\|_0^2.
\end{equation}
The unique solvability of the discrete problem follows (if $f\equiv 0$).

For the error estimation we introduce projections $\Pi: V\mapsto V^N$ and $P: W\mapsto W^N$. As usual,
instead of $u-u^N$ and $\bar q-\bar q^N$ we estimate $Pu-u^N$ and $\Pi\bar q-\bar q^N$, assuming that
we can estimate the projection errors. Subtraction
of the continuous and the discrete problem results in
\begin{subequations}
  \begin{align}
\eps(\nabla\cdot( \Pi\bar q-\bar q^N),w)+c(Pu-u^N,w)&=\eps(\nabla\cdot( \Pi\bar q-\bar q),w)+c(Pu-u,w),\\
\eps(\Pi\bar q-\bar q^N,\bar v)-\eps(\nabla\cdot \bar v,Pu-u^N)&=\eps(\Pi\bar q-\bar q,\bar v)-\eps(\nabla\cdot \bar v,Pu-u).
   \end{align}
\end{subequations}
Setting $\bar v:=\Pi\bar q-\bar q^N)=\bar \mu$ and $w:=Pu-u^N=\tau$ we obtain the error equation
\begin{equation}
\eps(\bar \mu,\bar \mu)+c(\tau,\tau)=\eps(\nabla\cdot( \Pi\bar q-\bar q),\tau)+c(Pu-u,\tau)+\eps(\Pi\bar q-\bar q,\bar \mu)-\eps(\nabla\cdot \bar \mu,Pu-u).
\end{equation}
From the error equation it is easy to derive a first order uniform convergence result in the energy norm (one could also
think about supercloseness similar as in \cite{LW00}). But we want to investigate, whether or not an estimate of the
type \eqref{est1} is possible.

If $P$ denotes the $L_2$ projection, we have
\[
(Pu-u,\tau)=0\quad {\rm and}\,\,(\nabla\cdot \bar \mu,Pu-u)=0,
\]
because $\nabla\cdot \bar \mu$ is piecewise constant for $\bar \mu\in V^N$. Therefore, from the right hand side of
the error equation two terms disappear and it follows
\begin{equation}\label{a}
\|\bar \mu\|^2_0 \preceq \eps\|\nabla\cdot( \Pi(\nabla u)-\nabla u)\|_0^2+\|\Pi(\nabla u)-\nabla u\|_0^2.
\end{equation}

Now let us denote by $\Pi^*$ the standard local projection operator into the Raviart-Thomas space
$V^N$. This operator satisfies
\begin{equation}\label{ort}
  (\nabla\cdot (\bar v-\Pi^*\bar v),w)=0 \qquad \text{for all } w\in W^N.
\end{equation}
Consequently, the choice $\Pi=\Pi^*$ would eliminate one more term in the error equation and thus in
\eqref{a}. But do we have for the projection error the desired estimate
\begin{equation}
\varepsilon^{1/4}\|\Pi^*(\nabla u)-\nabla u\|_0\preceq N^{-1}\ln N \quad ?
\end{equation}
The answer is no (see Lin and Stynes \cite{LS12}, page 2738). The reason lies in the fact that
$\Pi^*$ is applied to $\nabla u$ and its behavior near the transition point of the mesh is different
from the behavior of $u$ (a factor $\eps^{-1}$).

Therefore, Lin and Stynes define a modified interpolant $\Pi\bar v\in W^N$, such that
\begin{equation}
\varepsilon^{1/4}\|\Pi(\nabla u)-\nabla u\|_0\preceq N^{-1}\ln N
\end{equation}
(\cite{LS12}, Corollary 4.6). The operator $\Pi$ is defined differently for every component
of the solution decomposition. For the smooth part one takes simply $\Pi=\Pi^*$.

For the layer components, however, $\Pi^*$ is modified. Consider, for instance, the layer
component $w_1$ related to $\exp(-\sqrt{c}y/\sqrt{\eps})$. Then $\Pi$ and $\Pi^*$ differ only
in the small strip $R_1$ defined by
\[
R_1:=[0,1]\times[\lambda-h^*,\lambda]\qquad {\rm with}\quad \lambda=2\sqrt{\eps}\ln N/\sqrt{c}
\quad
{\rm and} \quad h^*=O(\sqrt{\eps}N^{-1}\ln N).
\]
On that strip we loose the property \eqref{ort}, therefore we have additionally to estimate
\begin{equation}
M_{1,R_1}:= \eps^{1/2}\|\nabla\cdot( \Pi(\nabla w_1)-\nabla w_1)\|_{0,R_1}.
\end{equation}
On $R_1$ we have $\|\Delta w_1\|\preceq \eps^{-1}N^{-2}$, consequently
\begin{equation}
\varepsilon^{1/2}\|\Delta w_1\|_{0,R_1}\preceq  \eps^{-1}N^{-2}\eps^{1/4}N^{-1/2}(\ln N)^{1/2}=
      \eps^{-1/4}N^{-5/2}(\ln N)^{1/2}.
\end{equation}
By construction the components of $\Pi ( \nabla w_1)$ satisfy $(\Pi (\nabla w_1))_1=0$ on $R_1$ and
$\|(\Pi (\nabla w_1))_2\|_\infty\preceq \eps^{-1}N^{-2}$. It follows
\begin{equation}
 \eps^{1/2}\|\nabla\cdot( \Pi\nabla w_1)\|_{0,R_1}\preceq \eps^{1/2}\frac{1}{h^*}\eps^{-1/2}N^{-2}(h^*)^{1/2}
       = \eps^{-1/4}N^{-3/2}(\ln N)^{-1/2}.
\end{equation}
Therefore
\begin{equation}
M_{1,R_1}\preceq\eps^{-1/4}N^{-3/2}.
\end{equation}
The other layer components of the solution decomposition of $u$ are treated similarly.
We obtain finally
\begin{equation}
\varepsilon^{1/4}\|\Pi\bar q-\bar q^N\|_0\preceq N^{-1}\ln N
\end{equation}
and
\begin{equation}
\varepsilon^{1/4}\|\nabla u-\bar q^N\|_0\preceq N^{-1}\ln N.
\end{equation}

\begin{rem}
It is well known \cite{AB85}, \cite{AC95} that mixed methods can be reformulated as
non-mixed formulations, more precisely as projected nonconforming methods. This allows
as well error estimates for certain nonconforming methods as the implementation of a
mixed method as nonconforming method.
\end{rem}


\begin{thebibliography}{AAA 00}

\bibitem{ALM15} Adler, J., MacLachlan, S., Madden, N.: A first-order system
Petrov-Galerkin discretisation for a reaction-diffusion problem on a fitted
mesh. to appear in IMA J. Num. Anal.


\bibitem{AC95} Arbogast, T., Chen, Z.:
   On the implementation of mixed methods as nonconforming methods
   for second order elliptic problems.
   Math. Comp., 64(1995), 943-972.

\bibitem{AB85} Arnold, D. N., Brezzi, F.:
   Mixed and nonconforming finite element methods.
   RAIRO, 19(1985), 7-32.

\bibitem{CT87} Crouzeix, M., Thomee, V.: The stability in $L_p$ and $W_p^1$ of
the $L_2$-projection onto
 finite element function spaces.
 Math. Comp., 48(1987), 521-532


\bibitem{DR97} Dobrowolski, M., Roos, H.-G.: A priori estimates for the
solution of convection-diffusion problems and interpolation on
Shishkin meshes. Z.~Anal.~Anwend., 16(1997), 1001-1012

\bibitem{FLR09} Franz, S.,Liu, F., Roos, H.-G., Stynes,M., Zhou, A.: The combination
technique for a two-dimensional convection-diffusion problem with exponential layers.
Appl. Math. 54(3), 2009, 203-223


\bibitem{FR14} Franz, S., Roos, H.-G.: Error estimates in a balanced norm for a
    convection-diffusion problme with two different boundary layers.
    Calcolo, 51(2014), 423-440

\bibitem{JF96} de Jager, E. M., Furu, J.:
   The theory of singular perturbations.
   North Holland 1996

\bibitem{KK05} Kellogg, R. B., Kopteva, N.: A singularly perturbed semilinear
reaction-diffusion problem in a polygonal domain.
J. diff. equ., 248(2010), 184-208


\bibitem{Li01} Li, J.:
   Uniform error estimates in the finite element method for a singularly
   perturbed reaction-diffusion problem.
    Appl. Numer. Math., 36(2001), 129-154.


\bibitem{Li97} Li, J.:
   Quasioptimal uniformly convergent finite element methods for the elliptic
   boundary layer problem.
    Computers Math. Applic., 33(1997), 11-22.

\bibitem{LW00} Li, J., Wheeler, M.F.: Uniform convergence and superconvergence of mixed
finite element methods on anisotropically refined grids.
SINUM, 38(2000), 770-798


\bibitem{LS12} Lin, R., Stynes, M.: A balanced finite element method for singularly
perturbed reaction-diffusion problems.
SINUM, 50(2012), 2729-2743

\bibitem{LMS09} Liu, F., Madden, N., Stynes, M., Zhou, A.: A two-scale sparse grid method for
singularly perturbed reaction-diffusion problems in two dimensions.
IMA J. Num. Anal., 29(4), 2009, 986-1007

\bibitem{MR15} Madden, N., Russell, S.: A multiscale sparse grid finite element method for
a two-dimensional singularly perturbed reaction-diffusion problem.
Adv. Comput. Math., 41(2015), 987-1014

\bibitem{MX15} Melenk, J. M., Xenophontos, C.: Robust exponential convergence of hp-FEM in
  balanced norms for singularly perturbed reaction-diffusion equations.
  Calcolo, 53(2016), 105-132

\bibitem{Os08} Oswald, P.: $L_\infty$-bounds for the $L_2$-projection onto
linear spline spaces. in: Recent advances in Harmonic Analysis and Applications,
Springer, New York 2013, 303-316


\bibitem{RST96} Roos, H.-G., Stynes, M., Tobiska, L.: Robust numerical methods
for singularly perturbed differential equations. \\Springer 2008

\bibitem{RS15} Roos, H.-G., Schopf, M.: Convergence and stability in balanced norms of finite element methods on Shishkin meshes
for reaction-diffusion problems. ZAMM, 95(6), 2015, 551-565

\end{thebibliography}
\end{document}